\documentclass[reqno]{amsart}

\usepackage{amssymb}

\newcommand\idd{\mathop {\fam 0 id}\nolimits}
\newcommand\Id{{\fam 0 I}}
\newcommand\Curr{\mathop {\fam 0 Cur}\nolimits}
\newcommand\Ann{\mathop {\fam 0 Ann}\nolimits}
\newcommand\torsion{\mathop {\fam 0 tor}\nolimits}

\newcommand\Cendl{\mathop {\fam 0 Cend}\nolimits^{\mathrm l}}
\newcommand\Cendr{\mathop {\fam 0 Cend}\nolimits^{\mathrm r}}
\newcommand\Cend{\mathop {\fam 0 Cend}\nolimits}
\newcommand\Hom{\mathop {\fam 0 Hom}\nolimits}
\newcommand\Choml{\mathop {\fam 0 Chom}\nolimits^{\mathrm l}}
\newcommand\Chomr{\mathop {\fam 0 Chom}\nolimits^{\mathrm r}}
\newcommand\gc{\mathop {\fam 0 gc}\nolimits}

\newcommand\Hmod{H\hbox{\rm -mod}}

\newcommand\Kerr{\mathop {\fam 0 Ker}\nolimits}
\newcommand\diag{\mathop {\fam 0 diag}\nolimits}

\newcommand\Ress{\mathop{ \fam 0 Res}}
\newcommand\GKdim{\mathop {\fam 0 GKdim}\nolimits}
\newcommand\rank{\mathop {\fam 0 rank}\nolimits}
\newcommand\End{\mathop {\fam 0 End}\nolimits}

\newcommand\length{\mathop{\fam 0 length}\nolimits}

\newcommand{\oo}[1]{\mathrel{{\circ }_{#1}} }

\newtheorem{thm}{Theorem}
\newtheorem{lem}[thm]{Lemma}
\newtheorem{prop}[thm]{Proposition}
\newtheorem{cor}[thm]{Corollary}

\theoremstyle{definition}
\newtheorem{defn}[thm]{Definition}
\newtheorem{exmp}[thm]{Example}
\newtheorem{rem}[thm]{Remark}

\numberwithin{equation}{section}
\numberwithin{thm}{section}

\begin{document}

\title[Simple associative conformal algebras of linear growth]{Simple associative conformal algebras \\ of linear growth}

\author{Pavel Kolesnikov}

\email{pavelsk@kias.re.kr}

\thanks{Partially supported by RFBR 05-01-00230 and SSc-2269.2003}

\address{Korea Institute for Advanced Study,
 207-43 Cheongnyangni 2-dong, Dongdaemun-gu, Seoul 130-722, Korea}

\keywords{conformal algebra, Weyl algebra, Gel'fand--Kirillov dimension}

\subjclass{16S99, 16P90, 16S32}

\begin{abstract}
We describe  simple  finitely  generated associative conformal algebras
of Gel'fand--Kirillov dimension one.
\end{abstract}

\maketitle

\section{Introduction}\label{Section1}

The theory of conformal algebras appeared as a formal language
describing the algebraic properties of the
operator product expansion (OPE) in two-dimensional conformal field theory
\cite{BPZ,FLM,K1,K3}.
In a few words, the origin of this formalism is as follows:
for any two local fields
$a(z), b(z)$ the commutator $[a(w),b(z)]$ could be written
as a (finite) distribution
with respect to derivatives of the delta-function:
\begin{equation}\label{OPE}
[a(w),b(z)]  = \sum\limits_{n\ge 0} c_n(z) \frac{1}{n!} \partial^n_z \delta(w-z),
\end{equation}
where $\delta(w-z) = \sum_{s\in {\mathbb Z}} w^{s} z^{-s-1}$.
The coefficients $c_n(z)$, $n\ge 0$,
of this distribution are considered as  new
``$n$-products" of $a(z)$ and $b(z)$. The algebraic
properties of these operations could be formalized by
a family of axioms, which gives the notion of a Lie conformal algebra~\cite{K1}.
Associative conformal algebras naturally come from representations
of Lie conformal algebras. Moreover, any of those Lie conformal
(super-) algebras appeared
in physics is embeddable into an associative one (this is not the case in general,
see \cite{Ro2}). So the investigation of associative conformal algebras
provides some information on the structure and representations of Lie
conformal algebras.

From the abstract point of view, a conformal algebra is a vector space $C$
endowed with a linear map $D:C\to C$ and with an infinite family of bilinear
operations $\oo{n}:C\times C \to C$
 ($n$ ranges through non-negative integers),
satisfying certain axioms.

In a more general context,
a conformal algebra is just an algebra in the pseudotensor
category \cite{BD} associated with the polynomial Hopf algebra $H= \Bbbk[D]$
(see also \cite{BDK}). This category consists of left $H$-modules endowed with
$H$-polylinear maps.
The notions of associativity,
commutativity, etc. are well-defined there.
If a conformal algebra $C$ is
a finitely generated $H$-module then $C$ is said to be
finite.

The structure theory of finite conformal algebras (and superalgebras)
was established in a series of works.
In \cite{CK2,DK,FK,FKR},
simple and semisimple finite Lie conformal (super)algebras were
described. Similar results for associative and Jordan conformal
algebras were obtained in  \cite{K3} and
\cite{Z1}, respectively.
In \cite{BDK}, the results of \cite{DK,K3} were generalized
for pseudotensor categories related to arbitrary cocommutative Hopf algebras
with finite-dimensional spaces of primitive elements.

There are two obvious ways how to move beyond the class of finite
conformal algebras. The first option is to consider free conformal
algebras and study their properties, as it was done in
\cite{Ro1,Ro2,Ro3}, see also \cite{BFK,BFK1,BFK2}. Another way is
to introduce a ``growth function" of a conformal algebra and move
into the next class with respect to the growth rate.
The present paper  develops the last
approach.

In
\cite{Re1}, the analogue of Gel'fand--Kirillov dimension ($\GKdim$)
for conformal algebras was proposed. As in the case of ordinary
algebras, finite conformal algebras have $\GKdim$~0, and there
are no conformal algebras with $\GKdim$ strictly between 0~and~1.
So the natural problem is to explore conformal algebras of
linear growth (i.e., those of  $\GKdim$~1). Since \cite{Re1},
this class of conformal algebras has been studied in
\cite{BKL1,BKL2,DK2,Re2,Z2}.

The difference between ordinary and conformal
algebras becomes apparent in the following
context. Although there are no (ordinary) finitely
generated simple associative algebras with
$\GKdim$~1 \cite{SW,SSW},  conformal algebras of this kind do exist.
There also exist infinite conformal algebras
with finite faithful irreducible representations.
It was conjectured in \cite{Z2} that these families of algebras coincide:

\smallskip
\noindent
{\bf Conjecture.}
An infinite associative conformal algebra
$C$ has a finite faithful irreducible representation
if and only if $C$ is a finitely generated simple
conformal algebra of $\GKdim$~1.
\smallskip

It was shown in \cite{Re1} that a finitely generated
simple associative conformal algebra
$C$ of linear growth
with a unit (i.e., with an element
$e\in C$ such that $e\oo{0}x=x$ for all $x\in C$ and
$e\oo{n}e =0$ for $n\ge 1$)
is isomorphic to the conformal algebra $\Cend_N$ of all
conformal endomorphisms \cite{DK,K3}
of the free $N$-generated $\Bbbk[D]$-module~$V_N$.
In this case, the conformal algebra $C$ has the
finite faithful irreducible
representation on~$V_N$.

In \cite{BKL1}, the converse statement was proved:
if a unital
associative conformal algebra $C$ has a finite faithful irreducible
representation,
then $C\simeq \Cend_N$. Therefore, the results of
\cite{BKL1,Re1} prove the Conjecture for unital conformal
algebras, but it is unclear how to join a unit to a conformal
algebra.

In \cite{Z2}, the Conjecture
(rather its ``if'' part) was confirmed for
conformal algebras with an idempotent
(i.e., $e\in C$ such that $e\oo{n}e =\delta_{n,0}e$, $n\ge 0$).

On the other hand, the classification of associative conformal algebras
with finite faithful irreducible representations
was proposed in \cite{K3},
partially constructed in \cite{BKL1},  and completed in \cite{Ko2}.
 In the present paper, we use the result of \cite{Ko2} in order
to prove the Conjecture in general, and get the explicit description of
finitely generated simple associative conformal algebras of linear growth.

\section{Associative conformal algebras
and their representations}\label{sec2}

From now on,
$\Bbbk $ is an algebraically closed field
of zero characteristic, $H$ is the polynomial algebra
$\Bbbk [D]$ with the ordinary derivation $\partial_D = \frac{d}{dD}$,
${\mathbb Z}_+$ denotes the set
of non-negative integers. We will also use the notation
$x^{(n)}$ for $x^n/n!$, $n\in {\mathbb Z}_+$.
It is convenient to set $x^{(n)}:=0$ for $n<0$.


\subsection{Conformal homomorphisms and products}\label{subsec2.1}

In this section, we introduce the basic definitions
from the theory of associative conformal algebras.
The exposition follows the one of \cite{BDK}, namely,
we adjust the theory of pseudoalgebras by making use
of commutativity of $H$. This point of view is also close to
\cite{DK,K3}, where the notion of a conformal algebra is introduced
via $\lambda$-brackets, but we avoid using ``external" variables.

Denote by $\Hmod $ the class of all unital left $H$-modules. Let
us consider the usual $D$-adic topology on $H$, i.e., suppose the
set of ideals $\{(D^n)\mid n\ge 0\}$ to be the system of basic
neighborhoods of zero in $H$. For any vector spaces $U$ and $V$
denote by $\Hom (U,V)$ the vector space of $\Bbbk $-linear maps
from $U$ to~$V$. Recall the notion of finite topology \cite{J} on
$\Hom(U,V)$: the family of basic neighborhoods of a map $\phi_0
\in \Hom (U,V)$ is presented by
\[
  \{\phi \in \Hom(U,V) \mid \phi(u_i)=\phi_0(u_i),\, i=1,\dots, m\},
\]
$u_1,\dots, u_m\in U$,
$m\ge 0$.
This topology turns $\Hom(U,V)$ into a topological vector space,
where a sequence $\{\phi_n\}_{n=0}^\infty \subset \Hom(U,V)$
converges to zero if and only if for any
$u_1,\dots, u_m\in U$,
$m\ge 1$, there exists
$N\ge 1$ such that $\phi_n(u_i)=0$, $i=1,\dots, m$, for all $n\ge N$.

\begin{defn}[\cite{DK}]\label{defn Chom}
Let $U,V\in \Hmod$. A continuous linear map
\begin{equation}\label{Chom}
a : H\to \Hom(U,V)
\end{equation}
is said to be a {\em left conformal homomorphism\/} from $U$ to $V$
if
\begin{equation}\label{sesqui-l}
[D, a(h)] = - a(\partial_D h), \quad h\in H.
\end{equation}
A {\em right conformal homomorphism\/} from $U$ to~$V$
is a continuous linear map (\ref{Chom}) satisfying
\begin{equation}\label{sesqui-r}
a(h)D  =-a(\partial_D h), \quad h\in H.
\end{equation}
\end{defn}

Denote by $\Choml(U,V)$ ($\Chomr(U,V)$) the vector space of all
left (right) conformal homomorphisms from $U$ to~$V$.
One may consider $\Choml(U,V)$ and $\Chomr(U,V)$ as $H$-modules
with respect to
\begin{eqnarray}
(Da)(h) = - a(\partial_D h), \quad a\in \Choml(U,V),\ h\in H,
  \label{l-Hmod} \\
(Da)(h) = D a(h) + a(\partial_D h),\quad a\in \Chomr(U,V),\ h\in H
 \label{r-Hmod}.
\end{eqnarray}

\begin{prop}[\cite{DK}]\label{propTorsion}
Let $u$ be an element of the $H$-torsion of $U$.
Then for any $a\in \Choml(U,V)$ or
$b \in \Chomr(U,V)$ we have $a(h)(u)=0=b(h)(u)$
for all $h\in H$.
\qed
\end{prop}

Given $U,V\in \Hmod$, the tensor product $U\otimes V$
could be considered an as $H$-module with respect to three
different structures. Let us denote these $H$-modules by
$U\otimes^{\mathrm o} V$ (outer), $U\otimes^{\mathrm l} V$
(left-justified), $U\otimes^{\mathrm r} V$ (right-justified):
\begin{eqnarray*}
D (u \otimes^{\mathrm o} v) &= & Du \otimes^{\mathrm o} v
  + u\otimes^{\mathrm o} Dv, \\
D (u\otimes^{\mathrm l} v) &=& Du \otimes^{\mathrm l} v, \\
D (u \otimes^{\mathrm r} v)  &=& u \otimes^{\mathrm r} Dv.
\end{eqnarray*}

The following definition is just a slightly modified one from
\cite{DK,K1,K3} (c.f. \cite{BDK}).

\begin{defn}\label{defn conf prod}
Let $U,V,W\in \Hmod$.
A map
\begin{equation}\label{conf prod}
\mu \in \Choml (U\otimes^{\mathrm r}V ,W),
\end{equation}
is said to be a $W$-valued
{\em conformal product\/} of $U$ and~$V$
if for any $h\in H$ the map $\mu(h)$
is a homomorphism of $H$-modules $U\otimes^{\mathrm o}V$ and~$W$.
\end{defn}

It is clear that any conformal product (\ref{conf prod})
could be considered as an element of $\Chomr(U\otimes^{\mathrm l} V,W)$,
so the notion of a conformal product is ``symmetric".
Denote the set of all $W$-valued conformal products of $U$ and $V$
by ${\mathrm P}(U,V; W)$.

One may interpret $\mu \in {\mathrm P}(U,V; W)$ as an $H$-linear map
$\mu_{\mathrm l}: U\to \Choml (V,W)$
or as
$\mu_{\mathrm r}: V\to \Chomr(U,W)$ via
\[
\mu(h)(u\otimes v) = [\mu_{\mathrm l}(u)(h)](v) =
 [\mu_{\mathrm r}(v)(h)](u),
 \quad h\in H,\ u\in U,\ v\in V.
\]

\begin{rem}\label{pseudo-remark}
Definition \ref{defn conf prod} provides a particular example
of an $H$-polylinear map in the corresponding
pseudotensor category \cite{BDK,BD}.
\end{rem}


\subsection{Conformal ${n}$-products}\label{subsec2.2}

Let us fix $U,V,W\in \Hmod$.
For any conformal homomorphisms
$a\in \Choml(V,W)$, $b\in \Chomr(U,W)$ and
for any vectors $u\in U$, $v\in V$, write
\begin{equation}\label{n-action}
a(D^{n})(v) = a\oo{n} v, \quad b(D^{n})(u)=u\oo{n} b,
\quad n\in {\mathbb Z}_+.
\end{equation}
The sequence of functions $a(D^{n})$ (resp., $b(D^{n})$),
$n\in {\mathbb Z}_+$, completely describes the
corresponding conformal homomorphism.

In particular, a $W$-valued conformal product~$\mu $ of $U$ and $V$
gives rise to the family of ordinary
linear operations ($n$-{\em products})
\begin{equation}\label{n-product}
\oo{n}: U\otimes V \to W, \quad n\in {\mathbb Z}_+,
\end{equation}
defined as
\begin{equation}
u\oo{n} v = \mu(D^n)(u\otimes v) = [\mu_{\mathrm l}(u)(D^{n})](v) =
 [\mu_{\mathrm r}(v)(D^{n})](u).
\end{equation}
It follows from Definition \ref{defn conf prod}
that
\begin{eqnarray}
& u\oo{n} v = 0, \quad \mbox{for $n$ sufficiently large},
                                             \label{C1}   \\
& u\oo{n} Dv = D(u\oo{n} v) + n u\oo{n-1} v,
                                             \label{C2}   \\
& Du\oo{n} v = -n u\oo{n-1} v.
                                             \label{C3}
\end{eqnarray}
Moreover, any sequence of $n$-products (\ref{n-product})
satisfying (\ref{C1})--(\ref{C3})
uniquely defines a conformal product in the sense
of Definition~\ref{defn conf prod}.

We will use the following notation:
for two subspaces $X\subseteq U$, $Y\subseteq V$
denote
\[
X\oo{\omega } Y
 = \sum\limits_{n\ge 0} (X\oo{n} Y) \subseteq W .
\]
If $Y$ is an $H$-submodule of $V$, then $X\oo{\omega} Y$ is also an
$H$-submodule of $W$ (see (\ref{C2}), (\ref{C3})).

Let us consider six modules $U_1,U_2,U_3,V_1,V_2,W\in \Hmod $
with conformal products
\begin{eqnarray}
\mu_{12} \in {\mathrm P}(U_1,U_2; V_1),
  \quad \nu_1 \in {\mathrm P}(V_1,U_3; W), \label{adj-1}\\
\mu_{23} \in {\mathrm P}(U_2,U_3; V_2),
  \quad \nu_2 \in {\mathrm P}(U_1,V_2; W).  \label{adj-2}
\end{eqnarray}
Each of these products is equivalent to a sequence of
$n$-products defined by (\ref{n-product}). We will denote these
operations uniformly by $\oo{n}$, $n\in \mathbb Z_+$.

One may consider the compositions of the conformal products
(\ref{adj-1}), (\ref{adj-2})
\[
\nu_1(\mu_{12}, \idd_{U_3}), \nu_2(\idd_{U_1}, \mu_{23}) :
H\otimes H \to \Hom (U_1\otimes U_2 \otimes U_3, W)
\]
defined as follows:
\begin{eqnarray}
\nu_1(\mu_{12}, \idd_{U_3})(h \otimes g) =
\nu_1(g)(\mu_{12}(h), \idd_{U_3}),
                             \label{eq (n)m} \\
\nu_2(\idd_{U_1}, \mu_{23})(h \otimes g) =
  \nu_2(h)(\idd_{U_1}, \mu_{23}(g)), \label{eq n(m)}
\end{eqnarray}
for any $h,g\in H$.

The maps (\ref{adj-1}), (\ref{adj-2}) may satisfy the
{\em associativity relation}. This relation has the following form:
\begin{equation}\label{conf prod ass}
\nu_1(\mu_{12}, \idd_{U_3}) = \nu_2(\idd_{U_1}, \mu_{23})\circ {\mathcal F},
\end{equation}
where
$\mathcal F$ is the formal Fourier transform \cite{BDK} given by
${\mathcal F} = \exp(-\partial_D \otimes D)$,
i.e.,
\[
\mathcal F(D^n\otimes D^m) = \sum\limits_{s\ge 0}
 (-1)^s\binom{n}{s} D^{n-s}\otimes D^{m+s}.
\]
Note that
${\mathcal F}^{-1} = \exp(\partial_D \otimes D)$,
so
\[
{\mathcal F}^{-1}(D^n\otimes D^m) = \sum\limits_{s\ge 0}
\binom{n}{s} D^{n-s}\otimes D^{m+s}.
\]

In terms of $n$-products, the relation (\ref{conf prod ass})
could be expressed as
\begin{equation}\label{conf-ass}
 (u\oo{n} v) \oo{m} w = \sum\limits_{s=0}^n
 (-1)^s \binom{n}{s} u\oo{n-s} (v\oo{m+s} w),
\end{equation}
or
\begin{equation}\label{conf-ass1}
u\oo{n} (v\oo{m} w) = \sum\limits_{s\ge 0} \binom{n}{s}
 (u\oo{n-s} v) \oo{m+s} w.
\end{equation}
for any  $n,m\in {\mathbb Z}_+$, $u\in U_1$,
$v\in U_2$, $w\in U_3$.
It is easy to note that the systems of relations
(\ref{conf-ass}) and (\ref{conf-ass1}) are equivalent.


\subsection{Conformal algebras and modules}\label{subsec2.3}

\begin{defn}[\cite{BDK}]\label{defn conf-alg}
An $H$-module $C\in \Hmod $ endowed with a conformal product
$\mu \in {\mathrm P}(C,C; C)$ is called a {\em conformal algebra}.
\end{defn}

The definition of a conformal algebra could be stated in terms of
$n$-products.

\begin{defn}[\cite{K1}]\label{defn conf-n}
A vector space $C$ endowed with a linear map $D$ and
with a family of linear maps
$\oo{n}: C\otimes C \to C$, $n\in {\mathbb Z}_+$,
is a conformal algebra if
(\ref{C1})--(\ref{C3}) hold.
\end{defn}

A conformal algebra $C$ is said to be associative if
$\mu(\mu, \idd_C) = \mu (\idd_C, \mu)\circ {\mathcal F}$,
i.e., if either of the relations (\ref{conf-ass}) or (\ref{conf-ass1})
holds for any $n,m\in {\mathbb Z}_+$.

Note that for any three subspaces $V_1$, $V_2$, $V_3$ of an associative
conformal algebra we have
\[
V_1\oo{\omega } (V_2\oo{\omega} V_3 )  = (V_1\oo{\omega } V_2)\oo{\omega} V_3  .
\]

A left (right) ideal $I$ of a conformal algebra $C$ is an
$H$-submodule $I\subseteq C$
such that $C\oo{\omega } I \subseteq I$
($I\oo{\omega } C\subseteq I$).
A conformal algebra $C$ is simple if $C\oo{\omega } C \ne 0$
and there are no
non-zero proper two-sided ideals of~$C$.
If $C$ is a finitely generated $H$-module, then $C$ is said to be
a finite conformal algebra.

\begin{rem}[\cite{DK}]\label{rem3.1}
It follows from Proposition~\ref{propTorsion} that any
simple conformal algebra is a torsion-free $H$-module.
\end{rem}

\begin{lem}\label{lem2.3.2}
If $C$ is a finitely generated associative conformal algebra
and%
\break
$\dim C/DC < \infty $,
then $C$ is finite.
\end{lem}

\begin{proof}
If $\dim C/DC = n< \infty$, then there exist
$a_1,\dots, a_n\in C$ such that any element $x\in C$
could be represented as
\[
x = x_0 + Dz_0, \quad x_0\in \Bbbk a_1 + \dots + \Bbbk a_n, \quad z_0\in C.
\]
Further, for any $s\ge 0$ there exist
 $x_s\in Ha_1 + \dots + Ha_n$, $z_s\in C$
such that
\begin{equation}\label{HSpan}
x = x_s + D^{s+1} z_s , \quad s\ge 0.
\end{equation}
Let $B$  be a finite set of generators
of $C$. Associativity relation implies that
\[
 C = \sum\limits_{b\in B}
    \left [
     Hb
      +
      \sum\limits_{t=0}^{N} H(C\oo{t}b)
    \right ]
\]
for sufficiently large~$N$. It follows from (\ref{HSpan})
that the finite set
$B\cup \{a_i\oo{t} b\mid i=1,\dots,n,\, b\in B, \, t=0,\dots, N\}$
generates $C$ as an $H$-module.
\end{proof}

Definitions \ref{defn conf-alg} and \ref{defn conf-n}
provide the axiomatic description of the following constructions
in formal distribution spaces.
Let $A$ be an algebra over $\Bbbk$ (non-associative, in general).
Then one may define the following operations on the space of
formal power series $A[[z,z^{-1}]]$:
\begin{equation}\label{Res}
a(z)\oo{n} b(z) = \Ress\nolimits_{w} a(w)b(z)(w-z)^n,\quad n\in {\mathbb Z}_+,
\end{equation}
where $\Ress_{w}F(w,z)$ means the coefficient at $w^{-1}$ of
a formal power series $F(w,z)$ in two variables.

Relations (\ref{C2}) and (\ref{C3}) hold with $D=d/dz$, but
the locality condition (\ref{C1}) does not hold in general. Two
series $a(z), b(z)\in A[[z,z^{-1}]]$ are called {\em local\/} if
$a(w)b(z)(w-z)^N = 0$ for some $N\in {\mathbb Z}_+$.
If a $D$-invariant subspace $C\subset A[[z,z^{-1}]]$ is closed under
all $n$-products (\ref{Res}) and consists of pairwise
mutually local series, then $C$ is a conformal algebra.
If $A$
is associative, then (\ref{conf-ass}) and (\ref{conf-ass1}) hold,
and any conformal algebra $C\subset A[[z,z^{-1}]]$ is associative.
The converse is also true (see, e.g., \cite{K3,Ro1}):
an (associative) conformal algebra $C$
in the sense of Definition~\ref{defn conf-n}
could be canonically embedded into the space of formal power series
$A[[z,z^{-1}]]$ over some (associative) algebra~$A$.

\begin{defn}[\cite{CK1,K1,K3}]\label{defn conf-mod}
Let $C$ be an associative conformal algebra with a
conformal product $\mu\in {\mathrm P}(C,C; C)$.
An $H$-module $V\in \Hmod $
endowed with a conformal product $\nu \in {\mathrm P}(C,V;V)$
is said to be a {\em left conformal $C$-module\/} if
$\nu (\mu, \idd_V) = \nu (\idd_C, \nu)\circ {\mathcal F}$.
Analogously, $V$ is said to be a {\em right conformal $C$-module\/}
if it is endowed with a conformal product
$\nu \in {\mathrm P}(V,C;V)$ such that
$\nu(\idd_V, \mu)\circ {\mathcal F} = \nu(\nu, \idd_C)$.
\end{defn}

It is clear how to express the last definition in terms of
$n$-products (see, e.g., \cite{CK1,DK}).

A (left) conformal $C$-module $V$ is called irreducible
if $C\oo{\omega }V \ne 0$ and there are no non-trivial
conformal $C$-submodules of~$V$.
If $V$ is a finitely generated $H$-module, then $V$ is said to
be a finite conformal $C$-module.


\subsection{Commutativity of conformal algebras}\label{subsec3.4}

For any $V,W\in \Hmod $, the $H$-modules
$\Choml(V,W)$ and $\Chomr(V,W)$ are isomorphic.
Namely,   for any
$a\in \Choml(V,W)$ the map $\tilde a : H\to \Hom (V,W)$ defined as
\begin{equation}\label{curly}
\tilde a(D^n) =
 \sum\limits_{s\ge 0} (-1)^{n+s} D^{(s)} a(D^{n+s}),
 \quad n\in \mathbb Z_+,
\end{equation}
belongs to $\Chomr(V,W)$ \cite{K3}.
The same relation
(\ref{curly}) defines the inverse map
$\Chomr(V,W) \to \Choml(V,W)$.
It is easy to check that
the map $\Choml(V,W)\to \Chomr(V,W)$, $a\mapsto \tilde a$,
is an isomorphism of $H$-modules.

Let $a\in \Choml(V,W)$, $b\in \Chomr(V,W)$
be conformal homomorphisms, $V,W\in \Hmod$.
For any $n\in {\mathbb Z}_+$,
denote
\[
\tilde a(D^n)(v)= \{a\oo{n} v\},
\quad
\tilde b(D^n)(v) = \{v\oo{n} b\}.
\]
By definition,
\begin{equation}\label{curly-n}
\{x\oo{n} y\} = \sum\limits_{s\ge 0}
 (-1)^{n+s} D^{(s)} (x\oo{n+s} y).
\end{equation}

Let $\mu \in {\mathrm P}(U,V; W)$ be a conformal product,
$U,V,W\in \Hmod$. Then the map
$\tilde \mu : H \to \Hom_{\Bbbk }(U\otimes V, W)$
constructed by the rule
\[
\tilde \mu(h)(u\otimes v) =
 [\tilde \mu_{\mathrm l}(h)(u)](v)
=[\tilde \mu_{\mathrm r}(h)(v)](u)
\]
could be considered as a $W$-valued conformal
product of $V$ and $U$: $\tilde \mu \in {\mathrm P}(V,U; W)$,
$\tilde \mu_{\mathrm l} = \mu_{\mathrm r}$,
$\tilde \mu_{\mathrm r} = \mu_{\mathrm l}$.
Using the notation introduced in \cite{BFKK},
denote the corresponding $n$-products by
$\{u\oo{n} v\}$, $u\in U$, $v\in V$, $n\in {\mathbb Z}_+$.
The relation between operations $\{\cdot \oo{n} \cdot \}$
and  $(\cdot \oo{n} \cdot)$ is given by (\ref{curly-n}).

\begin{prop}[\cite{BFKK,DK,K3}]\label{prop-L}
For any conformal products
(\ref{adj-1}), (\ref{adj-2})
such that the associativity relation (\ref{conf prod ass})
holds, we also have the following properties:
\begin{eqnarray}
u \oo{n} \{ v \oo{m} w\} & = &
 \{(u \oo{n} v) \oo{m} w\};   \label{eq2.2.1} \\
\{u \oo{n} (v \oo{m} w)\}  & = &
   \sum\limits_{s\ge 0} (-1)^s\binom{m}{s}
   \{\{u \oo{m-s} v\} \oo{n+s} w \}  ;
  \label{eq2.2.2} \\
\{ u \oo{n}\{v \oo{m} w \}\}  & = &
  \sum\limits_{s\ge 0}(-1)^s \binom{m}{s}
   \{\{ u\oo{n+s} v\} \oo{m-s} w \};
  \label{eq2.2.3} \\
\{u \oo{n} v\} \oo{m} w & = &
  \sum\limits_{s\ge 0} (-1)^s \binom{n}{s} u\oo{m+s}(v \oo{n-s} w),
  \label{eq2.2.4}
\end{eqnarray}
$u\in U_1$, $v\in U_2$, $w\in U_3$.
\qed
\end{prop}

A conformal algebra $C$ with a conformal product
$\mu $ is said to be commutative
if $\mu(h)(a\otimes b) = \tilde \mu(h)(b\otimes a)$
for all $h\in H$, $a,b\in C$,
i.e., if $a\oo{n} b = \{b\oo{n} a\}$.
Conformal algebra is commutative if and only if
it could be embedded into $A[[z,z^{-1}]]$ over
a commutative algebra~$A$~\cite{K3,Ro1}.

\begin{rem}[\cite{K3}]\label{l-rem}
Since a conformal product
$\mu \in {\mathrm P}(U,V; W)$,
$\mu : H \to \Hom (U\otimes V, W)$
is a continuous map, it is possible to define
$\mu(\exp(\alpha D))\in \Hom (U\otimes V, W)$, $\alpha \in \Bbbk$.
This operation is denoted by
$(\cdot _{\alpha} \cdot): U\otimes V \to W$ \cite{DK,K3}:
\[
(u {}_\alpha v) = \sum\limits_{n\ge 0} \alpha^{(n)} (u\oo{n}v).
\]
Analogously, we denote
\[
\{ u {}_\alpha v \}= \tilde \mu(\exp(\alpha D))(u\otimes v)
 = \sum\limits_{n\ge 0} \alpha^{(n)} \{u\oo{n}v\}.
\]
\end{rem}

The following proposition provides a useful description
of irreducible modules over associative conformal algebras.

\begin{prop}[\cite{BKL1}]\label{prop2.1}
Let $C$ be an associative conformal algebra, and let $V$ be an irreducible
(left) $C$-module. Then there exist $u\in V$, $\alpha \in \Bbbk $
such that $\{C_{\alpha } u\}=V$.
\qed
\end{prop}


\subsection{Algebra of $0$-multiplications}

Let $C$ be an associative conformal algebra with a
conformal product~$\mu$.
Then  for every $a\in C$ the operator
$a(0)=\mu_{\mathrm l}(1)(a) \in \End C $, $a(0) :  x\mapsto a\oo{0}x$,
$x\in C$,
is a homomorphism of right
$C$-modules (see (\ref{C3}), (\ref{conf-ass1})).
The set
\begin{equation}\label{eq2.1-A0}
A_0 = A_0(C)= \{a(0)\mid a\in C \}
\end{equation}
is an associative subalgebra of $\End C$:
$a(0)b(0) = (a\oo{0} b)(0)$.

\begin{lem}\label{lem2.3.1}
If $C$ is a finitely generated conformal algebra
and
\[
a_1\oo{0} C + \dots + a_m\oo{0} C = C
\]
for some
$a_1,\dots, a_m \in C$,
then $A_0$ is a finitely generated algebra.
\end{lem}

\begin{proof}
Let $B\subset C$ be a finite set of generators.
Fix an upper bound of locality on $B$, i.e., a number $N$
such that
$a\oo{n} b = 0$ for any
 $n> N$, $a,b,\in B$.
Denote by $I$ the finite set of symbols
$\{1,\dots,m\}$, and let $I^*$ be
the set of finite words in $I$
(by $\varepsilon $ we denote the empty word).

For any $b\in B$ we may construct a family
of elements $\{b_{(w)}\in C \mid w\in I^*\}$,
such that
$b_{(\varepsilon)}=b$
and
for any $u\in I^*$ we have
\[
b_{(u)} = a_1\oo{0}b_{(u1)} + \dots +a_m\oo{0}b_{(um)}.
\]
Such a family is not unique, but we may fix any of them.

For all $l,r,n\in {\mathbb Z}_+$, denote by
$B_n^{(l,r)}$
the set of all elements  of the form
\[
(a_{i_1}\oo{s_1} \dots \oo{s_{l'-1}}a_{l'}
   \oo{s_{l'}}  b_{(u)} \oo{p_1}
  a_{j_1}\oo{p_2} \dots \oo{p_{r'} } a_{j_{r'} }),
\]
where $i_k,j_t=1,\dots, m$, $s_k, p_t=0,\dots, N$,
($k=1,\dots, l'$, $t=1,\dots, r'$),
$b\in B$, $u\in I^*$, $\length(u)\le n$,
$l'\le l$, $r'\le r$, with
all possible bracketing schemes.
It is clear that $\bigl| B_n^{(l,r)}\bigr|<\infty $
and $B=B_0^{(0,0)}\subseteq B_n^{(l,r)}$
for all $l,r,n\in {\mathbb Z}_+$.

First, let us show that
\begin{equation}\label{local claim}
 B_{n_1}^{(l_1,r_1)}\oo{n}  B_{n_2}^{(0,r_2)}
 \subseteq
   B_{n_1}^{(l_1, r_1+ n)} \oo{0}
   B_{n+n_2}^{(n,r_2)},
 \quad n\le N.
\end{equation}
For $n=0$ the relation (\ref{local claim}) is trivial, so proceed
by induction on~$n$. Let $x\in   B_{n_2}^{(0,r_2)}$,
$y\in   B_{n_1}^{(l_1,r_1)}$.
By definition (see also (\ref{conf-ass})),
there exist $x_1,\dots, x_m \in   B_{n_2+1}^{(0,r_2)}$
such that
$x=a_1\oo{0}x_1 + \dots + a_m\oo{0} x_m$.
Then
\begin{eqnarray}
y   \oo{n+1}  x  & = & \sum\limits_{j=1}^m y\oo{n+1} (a_j \oo{0} x_j)
   \nonumber \\
  & = &
 \sum\limits_{j=1}^m
  \left[
    \sum\limits_{s=0}^n \binom{n+1}{s} (y\oo{n+1-s} a_j) \oo{s} x_j \right]
 \nonumber \\
   &  + &
 \sum\limits_{j=1}^m    (y\oo{0}a_j)\oo{n+1} x_j   .
  \label{local claim1}
\end{eqnarray}
Since $y\oo{n+1-s} a_j \in   B_{n_1}^{(l, r_1+1)}$,
the elements
$(y\oo{n+1-s} a_j)\oo{s} x_j$ lie in
the space
$  B_{n_1}^{(l_1, r_1+n+1)}\oo{0}   B_{n_2+n+1}^{(n,r_2)}$
for all $s=0,\dots, n$.
The last summand of (\ref{local claim1}) clearly lies in
$  B_{n_1}^{(l_1,r_1)}\oo{0}   B_{n_2+n+1}^{(1,r_2)}$.

It follows from (\ref{conf-ass1}) that any element of $C$ could be
presented as an $H$-linear combination
of left-normed words in $B$, i.e., conformal monomials
of the form
\[
(\dots (b_1\oo{n_1} b_2)\oo{n_2} \dots \oo{n_k} b_{k+1}),
\quad b_i\in B=B_0^{(0,0)}, \quad n_i\le N,
\]
generate $C$ as an $H$-module.
Relation (\ref{local claim}) implies that such a monomial lies
in
\[
\bigl(  B_{0}^{(0,n_1)} \oo{0}   B_{n_1}^{(n_1,n_2)} \oo{0}
 \dots \oo{0}   B_{n_{k-1}}^{(n_{k-1},n_k)} \oo{0}   B_{n_k}^{(n_k,0)} \bigr),
\]
hence, the finite set
$\bigl\{ b(0) \mid b\in   B_N^{(N,N)} \bigr\}$
generates associative algebra~$A_0$.
\end{proof}

\section{Gel'fand--Kirillov dimension of conformal algebras}\label{sec4}

Let $C$ be a finitely generated conformal algebra (not necessarily
associative), and let $\{a_i \mid i\in S\} \subset C$ be a finite
system of generators ($|S|<\infty $). Denote by $C_n$, $n\ge 1$,
the $H$-linear span of all conformal monomials of the form
$(a_{i_1}\oo{m_1} \dots \oo{m_{k-1}} a_{i_k})$, $k=1,\dots, n$,
with all possible bracketing schemes. For a fixed bracketing
scheme, there exist only a finite number of non-zero monomials, so
$C_n$ is finitely generated as an $H$-module. In particular, the
number $d_n =\rank C_n $ is finite. It is also clear that
$C=\bigcup_{n\ge 1} C_n$.

From now on, we consider associative conformal
algebras only. Throughout the rest of the paper,
the term ``conformal algebra"  means ``associative conformal algebra".

Let $V$ be a left $C$-module over a
finitely generated conformal algebra $C$. Assume that
$V$ is a finitely generated  $C$-module, and choose a finite system of generators
$\{v_j\mid j\in T\}$, $|T|<\infty $, of $V$ over $C$.
 Denote by $V_1$ the $H$-linear span of
 $\{v_j \mid j\in T\}$. For $n>1$, define
\begin{equation}\label{3.1*}
V_n =   V_{n-1} + C_1\oo{\omega }V_{n-1} = V_1 + C_{n-1}\oo{\omega } V_1.
\end{equation}
It is clear that
$V_n$ is a finitely generated $H$-module ($\rank V_n <\infty $),
and $V=\bigcup_{n\ge 1}V_n$.

The following definition was originally introduced for the regular module,
i.e., for a (non-associative, in general) conformal algebra itself.

\begin{defn}[\cite{Re1}]\label{GKdim}
Let $C$ be a finitely generated conformal  algebra,
and let $V$ be a finitely generated left $C$-module.
Then the value
\[
\GKdim_C V  =
\limsup\limits_{n\to \infty} \log_n (\rank V_n)\in
 {\mathbb R}_+ \cup \{\infty \}
\]
is called the {\em Gel'fand--Kirillov dimension\/} of $C$-module~$V$.

If $V=C$ is the regular module, then
$\GKdim C = \GKdim_C C$  is called the Gel'fand--Kirillov dimension of~$C$.
\end{defn}

The similar construction is valid for right $C$-modules.

The definition of Gel'fand--Kirillov dimension could be expanded
to  infinitely generated conformal algebras
and modules in the usual way:
\[
\GKdim _C V =\sup\limits_{C'\subseteq_{f.g.} C, \, V'\subseteq_{f.g.} V} \GKdim _{C'} V'.
\]

As in the case of usual algebras, Gel'fand--Kirillov dimension does not depend
on the system of generators of $C$ (and of $V$). If $V'$ is a conformal
 submodule or homomorphic image of $V$,
 then $\GKdim_C V' \le \GKdim_C V$ (see \cite{Re1}).
It is also easy to obtain the following properties
(well-known for usual algebras).

\begin{lem}\label{lem4.1}
Let $C$ be  a finitely generated conformal algebra,
and let $V$ be a finitely generated $C$-module.
Then

{\rm (i)} $\GKdim_C V = 0$ if and only if $V$ is  a finite
conformal module;

{\rm (ii)} $\GKdim_C V>0$ implies $\GKdim_C V \ge 1$.
\end{lem}

\begin{proof}
Although the statements are clear, let us state the proof in the case of
right modules. First, consider the case when $V$ is a torsion-free
$H$-module, i.e., $\torsion V =0$.
Assume that $\rank V_n = \rank V_{n+1}$ for some $n\ge 1$.
Then there exists $f\in H$ such that $fV_{n+1 } \subseteq V_n$,
so the $H$-module $V_{n+1}/V_n$ coincides with its torsion.
Since every generator $a\in C_1$ defines a conformal
homomorphism in $\Chomr(V_{n+1}/V_n, V_{n+2}/V_{n+1})$
by the rule
\[
( v+V_{n})\oo{m} a  = v\oo{m}a + V_{n+1}, \quad v\in V_{n+1},\ m\in {\mathbb Z}_+,
\]
Proposition \ref{propTorsion}
implies $V_{n+1}\oo{\omega } C_1 \subseteq V_{n+1}$.
Therefore, the ``right analog" of (\ref{3.1*})
implies $V_{n+1}=V_{n+2}=\dots = V$, so $\rank V<\infty$.

We have proved that either $\rank V<\infty $ or
the sequence
$\rank V_n$, $n\ge 1$, is strictly increasing.
In the last case,
$\rank V_n \ge n$ for every $n\ge 1$,
so $\GKdim _C V \ge 1$.

If $\torsion V \ne 0$, then we may apply the arguments above to the
torsion-free module $\overline V = V/\torsion V$. Since
$\GKdim_C \overline V \le \GKdim_C V$, it is sufficient to show
that $\rank (\overline V)<\infty$ implies $\rank V <\infty$.
Indeed, if $\overline V$ is finitely generated as an $H$-module,
then there exist $v_1,\dots, v_m \in V$ such that an arbitrary
$v\in V$ could be presented as
\begin{equation}                              \label{eq-3_pr}
v = f_1v_1 + \dots + f_m v_m + u,
 \quad  f_i\in H, \ u\in \torsion V.
\end{equation}
Since $V = \bigcup_{n\ge 1} V_n$,
there exists $n\ge 1$ such that $v_1,\dots,v_m \in V_n$.
In particular, every element $v\in V_{n+1}$ can be presented as (\ref{eq-3_pr}),
so by Proposition~\ref{propTorsion} we have
$V_{n+1}\oo{\omega} C_1 \subseteq V_{n+1}$.
Hence, $\rank V <\infty $ as above.
\end{proof}

\begin{cor}\label{corGK}
{\rm (i)}
Let $V$ be an arbitrary conformal module over a
conformal algebra $C$.  Then $\GKdim _C V<1$ implies $\GKdim_C V =0$.

{\rm (ii)}
There are no conformal algebras of Gel'fand--Kirillov dimension strictly between
one and zero.
\qed
\end{cor}

It is well-known for usual algebras, that if
an ideal $I$ of an associative algebra $A$ contains
a regular element, then $\GKdim A/I \le \GKdim A - 1$.
We will use an analogous fact for conformal modules.

Let us consider a finitely generated conformal algebra $C$,
and fix a system of generators.
A homomorphism of (right) $C$-modules
$L: C\to C$ is called {\em bounded\/} if there exists $m=m(L)\ge 0$ such
that
$L(C_n) \subseteq C_{n+m}$ for any $n\ge 0$.
The set of all bounded homomorphisms is an associative subalgebra
${\mathcal L}(C)$
of
$\End C$.
For example, $H$-linear maps of the form
$x\mapsto \alpha x + a\oo{0}x$, $\alpha \in \Bbbk $,
$a\in C$,
lie in ${\mathcal L}(C)$.

For every $L\in {\mathcal L}(C)$ consider
\begin{equation}\label{im+ker}
I_L = L(C) + \bigcup\limits_{n\ge 1}\Kerr_C L^n \subseteq  C.
\end{equation}
It is clear that $I_L$ is a right ideal of $C$.

\begin{prop}\label{MainProp}
Let $C$ be a
finitely generated
conformal algebra such that $C \oo{\omega} C= C $, and let $L\in {\mathcal L}(C)$.
If $I_L\ne C$, then
there exists an irreducible finitely generated right $C$-module $V$
such that
\[
\GKdim _C V \le \GKdim C -1.
\]
\end{prop}

\begin{proof}
Let
$J_1 \subseteq J_2 \subseteq J_3 \subseteq \dots $
be an increasing chain of proper right ideals of $C$.
If $J_\omega = \bigcup_{k\ge 1} J_k $ is equal to $C$,
then there exists $m\ge 1$ such that
$J_m $  contains all generators of $C$. In this case $J_m = C$
is not a proper ideal. Therefore, $J_\omega $ is a proper right ideal of $C$.
By the Zorn lemma, we conclude that every proper right ideal
of $C$ can be embedded into a maximal right ideal of $C$.
Let $I$ be a maximal right ideal which contains $I_L \ne C$.

Note that $V = C/I$ is an irreducible finitely generated right $C$-module
($V\oo{\omega }C \ne 0$ since $I\not\supseteq C\oo{\omega } C = C$,
this is the only place where we we use this condition).
In particular, $V$ is a torsion-free $H$-module.
Let us fix a finite system of generators $\{a_i \mid i\in S\}$
for $C$,
and consider $V$ to be generated by the set
$\{\bar a_i=a_i+I \mid i\in S\}$.

Let $d(n) = \rank V_n$,
$n\ge 1$, and let
$\bar u_1, \dots ,\bar u_{d(n)}$ be an $H$-basis of $V_n$.
Then the set of elements
\[
\{ u_j, L(u_j) ,\dots, L^n(u_j) \mid j=1,\dots, d(n) \} \subset C_{n+mn},
\quad m=m(L),
\]
is $H$-linearly independent.
Indeed,
if there exist some polynomials
$f_{k,j}$,
$k=0,\dots, n$,
$j=1,\dots, d(n)$,
such that
\begin{equation}\label{sum-C}
\sum\limits_{j} f_{0,j} u_j +
 \sum\limits_{j} f_{1,j} L(u_j)
 + \dots +
 \sum\limits_{j} f_{n,j} L^n(u_j) =0,
\end{equation}
then
$\sum_j f_{0,j}u_j \in L(C) \subseteq I$,
so $f_{0,j}=0$ for all $j=1,\dots, d(n)$.
The relation (\ref{sum-C}) implies that
\[
\sum\limits_{j} f_{1,j} u_j
 + \dots +
 \sum\limits_{j} f_{n,j} L^{n-1}(u_j) \in \Kerr_C L \subseteq I,
\]
so
$\sum_{j} f_{1,j} u_j \in I+L(C) =I$,
and all $f_{1,j}$ are zero.
In the same way,  $f_{k,j}=0$
for all
$k=0,\dots, n$,
$j=1,\dots, d(n)$.

Since $C_{n+mn}$ contains at least $(n+1)d(n)$ $H$-linearly independent
elements, we may conclude that
$\rank C_{n+mn} \ge n d(n)$.
Now it is left to apply the usual arguments to deduce
$\GKdim C \ge 1+ \GKdim_C V$.
\end{proof}

\begin{rem}\label{rem4.1}
The same statement could be proved for ``right'' bounded operators,
i.e., for homomorphisms $R: C\to C$ of left $C$-modules,
such that $R(C_n)\subseteq C_{n+m}$ for some $m\ge 0$ and for all $n\ge 1$.
For example, the maps
$x\mapsto \alpha x + \{x\oo{0}a\}$,
$\alpha \in \Bbbk$,
$a\in C$,
satisfy these conditions.
\end{rem}

\begin{lem}\label{lem5.1}
Let $C$ be a conformal algebra, and let $A_0 = A_0(C)\subseteq {\mathcal L}(C)$
be the algebra of $0$-multiplications defined by
(\ref{eq2.1-A0}). If $C$ is a torsion-free $H$-module, then
$\GKdim A_0\le \GKdim C$.
\end{lem}

\begin{proof}
Consider a subalgebra
$A_0' \subseteq A_0$ generated by a finite family of elements
$a_1(0), \dots , a_n(0)\in A_0$, $n\ge 1$, and let
$C'$ be the conformal subalgebra of $C$ generated by
$\{a_1, \dots , a_n\}$.
It is clear that $A_0' \subseteq A_0(C')$.

It is sufficient to show that if $u_1(0), \dots , u_d(0) \in A_0(C')$,
$d\ge 1$, are
$\Bbbk $-linearly independent, then $u_1,\dots, u_d\in C'$
are $H$-linearly independent. Assume
\[
f_1(D) u_1 + \dots + f_d(D)u_d = 0
\]
for some $f_1,\dots, f_d \in H$. Then
\[
f_1(0)u_1(0) + \dots + f_d(0)u_d(0) = 0.
\]
Since $C'$ is a torsion-free $H$-module, there exists $i\in\{1,\dots, d\}$
such that $f_i(0)\ne 0$, so $u_1(0),\dots, u_d(0)$ are linearly
dependent over $\Bbbk$.

Therefore, $\GKdim A_0' \le \GKdim C' \le \GKdim C$ for any
finitely generated subalgebra $A_0'$ of $A_0$. Thus, $\GKdim A_0\le \GKdim C$.
\end{proof}

\section{The conformal algebra $\Cend_N$ and its irreducible
 subalgebras}\label{Section3}

Let $V$ be a unital left $H$-module.
Denote the $H$-modules
$\Choml (V,V)$ and $\Chomr(V,V)$
by $\Cendl V$ and $\Cendr V$, respectively.
The vector spaces $\Cendl V$, $\Cendr V$
(of left and right {\em conformal
endomorphisms})
are considered as  $H$-modules
with respect to (\ref{l-Hmod}), (\ref{r-Hmod}).
There exist natural conformal products
\[
\nu_1 \in {\mathrm P} (\Cendl V, V; V),
\quad
\nu_2 \in {\mathrm P} (V,\Cendr V; V)
\]
defined as follows:
\[
\nu_1(h)(a\otimes v) = a(h)(v),
\quad
\nu_2(h)(v\otimes b) = b(h)(v),
\]
$h\in H$,
$v\in V$,
$a\in \Cendl V$,
$b\in \Cendr V$.

If $V$ is a finitely generated $H$-module, then $\Cendl V$ and $\Cendr V$
can be endowed with conformal products $\mu_1$ and $\mu_2$,
respectively, such that
\[
\nu_1(\mu_1, \idd_V) = \nu_1 (\idd_{\Cendl V}, \nu_1)\circ {\mathcal F},
\quad
\nu_2(\idd_V, \mu_2) = \nu_1 (\nu_2, \idd_{\Cendr V})\circ {\mathcal F}^{-1}.
\]
Therefore, $\Cendl V$ and $\Cendr V$ are
(associative) conformal algebras (see, e.g., \cite{BDK,DK,K3}).

Let $V_N$ be a free $N$-generated $H$-module, $N\ge 1$.
The conformal algebra $\Cendl V_N = \Cendl_N$ can be presented
as follows.
For a fixed $H$-basis $\{e_1,\dots, e_N\}$
of $V_N$, one may define the operation $\partial_D$ on $V_N$
as on $H\otimes \Bbbk^N$ ($\partial_D e_i = 0$, $i=1,\dots, N$).
For any matrix $A=A(x)\in M_N(\Bbbk [x])$ with polynomial
entries define the conformal endomorphism
\begin{equation}\label{Cend}
A : H \to \End V_N,
\quad
A(x)(D^n): u \mapsto A(D)\partial_D^n(u),
\end{equation}
where $n\in {\mathbb Z}_+$.

\begin{prop}[\cite{DK,K1,K3,Re1}]\label{prop CendN}
The $H$-module generated by the conformal endomorphisms
(\ref{Cend}), $A\in M_N(\Bbbk[x])$,
is isomorphic to $\Cendl_N$.
Therefore, the conformal algebra $\Cendl_N$ could be
identified with $M_N(\Bbbk[D,x])\simeq H\otimes M_N(\Bbbk[x])$,
where the family of $n$-products is defined by
\begin{equation}\label{n-prod-Cend}
A(x)\oo{n} B(x) = A\partial^n_x(B(x)), \quad n\in {\mathbb Z}_+.
\end{equation}
\end{prop}

The map $a\mapsto \tilde a$ given by (\ref{curly}) provides
an anti-isomorphism of conformal algebras $\Cendl_N$ and $\Cendr_N$.
Namely, we have
\begin{equation}\label{anti-iso}
\tilde a \oo{n} \tilde b = \widetilde{\{b\oo{n} a\}},
\quad a,b\in \Cendl_N, \ n\in {\mathbb Z}_+.
\end{equation}
Combining this map with an arbitrary anti-automorphism
of $\Cendl_N$
(see \cite{BKL1} or \cite{Ko2}),
we get an isomorphism of conformal algebras
$\Cendl_N$  and $\Cendr_N$. From now on,
we will denote this conformal algebra by $\Cend_N$.

It is clear from (\ref{Cend}) that
\begin{equation}\label{Weyl-operator}
\{a(h)\mid a\in \Cend_N,\, h\in H\} = M_N(W),
\end{equation}
where $W$ is the first Weyl algebra (see \cite{Ko2} for details):
\[
W = \Bbbk \langle D, \partial_D \mid \partial_D D - D\partial_D =1 \rangle
\subset \End V_N.
\]

\begin{exmp}\label{exmp CurN}
Consider
$S_0 = M_N(\Bbbk[\partial_D]) \subset M_N(W)$.
The set
\begin{equation}\label{Curr eq}
\{a\in \Cend_N \mid a(H)\subset S_0\} = \Curr_N
\end{equation}
is a conformal subalgebra of $\Cend_N$ called
the {\em current\/} subalgebra. The image of $\Curr_N$ in $M_N(\Bbbk[D,x])\simeq \Cend_N$
(see Proposition~\ref{prop CendN}) is given by $M_N(\Bbbk[D])$.
\end{exmp}

\begin{exmp}\label{exmp CendNQ}
For any matrix $Q=Q(D)\in M_N(\Bbbk[D])$
the set
\begin{equation}\label{CendNQ}
\{a\in \Cend_N \mid a(H)\subset M_N(W)Q \} = \Cend_{N,Q}
\end{equation}
is a conformal subalgebra of $\Cend_{N}$. This conformal subalgebra
could be identified with $M_N(\Bbbk[D,x])Q(-D+x)$.
\end{exmp}

\begin{rem}[\cite{BKL1,K3,Z1}]\label{remark4_GK}
Conformal algebras $\Cend_{N,Q}$,
$N\ge 1$, $\det Q\ne 0$, are simple, finitely generated, and
$\GKdim \Cend_{N,Q} = 1$. Moreover,
these are infinite conformal algebras with finite
faithful irreducible modules.
\end{rem}

The following statement was conjectured in \cite{K3}.
In \cite{Ko2}, it was proved for left modules, but it also
holds in the case of right modules since $\Cendl_N \simeq \Cendr_N$.

\begin{thm}[\cite{Ko2}]\label{thmFFF}
Let $C$ be a conformal algebra with a finite faithful
irreducible module.
Then either $C\simeq \Curr_N$ or $C\simeq \Cend_{N,Q}$,
for some $N\ge 1$, where
$Q$ has the canonical diagonal form, i.e.,
$Q=\diag(f_1,\dots, f_N)$,
all $f_i$ are non-zero monic polynomials, and
$f_1\vert f_2\vert \ldots \vert f_N$.
\qed
\end{thm}

\begin{thm}[\cite{Ko2}]\label{thmFF}
Let $C$ be a simple associative conformal algebra
with a non-trivial finite module.
Then either $C$ is finite, or
$C\simeq\Cend_{N,Q}$ as in Theorem~\ref{thmFFF}.
\qed
\end{thm}

\begin{cor}\label{cor4.1}
Let $C$ be a finitely generated simple conformal algebra
such that $\GKdim C =1$.
Then either
$C\simeq \Cend_{N,Q}$
($N\ge 1$, $\det Q\ne 0$)
or
for any $L\in {\mathcal L}(C)$ there exists $n\ge 1$
such that
$L(C) + \Kerr_C L^n =C$.
\end{cor}

\begin{proof}
Suppose there exists $L\in {\mathcal L}(C)$
such that for any $n\ge 1$ the ideal $L(C)+\Kerr_C L^n\subset C$ is proper
(see the proof of Proposition~\ref{MainProp}).
Since $C$ is a finitely generated conformal algebra, the union
of these ideals (equal to $I_L$) is also proper.
Then by Proposition~\ref{MainProp}
we conclude that there exists an irreducible
finitely generated $C$-module $V$
such that $\GKdim_C V =0$.
By Lemma \ref{lem4.1}, $V$ is finite,
so Theorem \ref{thmFF} implies that $C$
is isomorphic to
$\Cend_{N,Q}$, $N\ge 1$, $\det Q\ne 0$.
\end{proof}

We are going to show that the second case is impossible,
so $\Cend_{N,Q}$, $N\ge 1$, $\det Q\ne 0$, exhaust
all simple finitely generated associative conformal algebras
of linear growth.

\section{Classification theorem}\label{Section5}

The purpose of this work is to prove the following statement.

\begin{thm}\label{thm5.1}
Let $C$ be a finitely generated simple associative conformal
algebra such that $\GKdim C =1$. Then $C\simeq \Cend_{N,Q}$,
$N\ge 1$, $\det Q \ne 0$.
\end{thm}

Throughout this section,
$C$ is assumed to satisfy the conditions of Theorem~\ref{thm5.1}.
This theorem proves the conjecture from
\cite{Z2} (see Section~\ref{Section1}) and generalizes the following results.

\begin{thm}[\cite{Re1}]\label{thm5.2}
If there exists an element
$e\in C$ such that
$e\oo{0}a=a$ for any $a\in C$, and
$e\oo{n}e =0$ for any $n\ge 1$,
then $C\simeq \Cend_N$, $N\ge 1$.
\end{thm}

Such an element $e\in C$ is called a {\em unit\/} of a conformal algebra.
Note that a conformal unit is not unique.

\begin{thm}[\cite{Z2}]\label{thm5.3}
If $C$ contains an element~$e$ such that
$e\oo{n}e=\delta_{n,0}e$
for any $n\ge 0$,
then $C\simeq \Cend_{N,Q}$, $N\ge 1$,
$Q=\diag (1, f_2,\dots, f_N)$.
\end{thm}

An element $e\in C$ satisfying the condition described in the last statement
is called an {\em idempotent\/} of a conformal algebra.

We will deduce a statement which is more general than Theorem \ref{thm5.2}.
Moreover, we will prove Theorem~\ref{thm5.3} using arguments different
from~\cite{Z2}.

Let us first sketch the idea of the proof of Theorem~\ref{thm5.1}.
Recall the
associative algebra
\[
A_0= A_0(C)=\{a(0): x\mapsto a\oo{0}x \mid  a\in C\} \subseteq {\mathcal L}(C)
\]
which was defined in
(\ref{eq2.1-A0}).
By Lemma~\ref{lem5.1} $\GKdim A_0 \le 1$.

If there exists $L\in {\mathcal L}(C)$ such that
$L(C)  + \Kerr_C L^n \ne C$
for all $n\ge 1$, then we may use Corollary~\ref{cor4.1}.

Therefore, it is sufficient to consider
 the case when
for any $L\in {\mathcal L}(C)$ there exists $n\ge 1$
such that
$L(C)  + \Kerr_C L^n = C$.
In particular, for any $a\in C$ we may assume
\begin{equation} \label{eq5.0}
a\oo{0} C  + \Kerr_C a(0)^n = C
\end{equation}
for an appropriate $n\ge 1$.
Relation (\ref{eq5.0}) implies
\begin{equation}\label{eq5.1}
a(0)A_0 + \Ann_{A_0} a(0)^n = A_0.
\end{equation}
It remains to show that either of (\ref{eq5.0}) or (\ref{eq5.1})
leads to a contradiction.

\begin{prop}\label{prop5.2}
If $A$ is a finitely generated algebra
of at most linear growth (i.e., $\GKdim A \le 1$),
 and for any $a\in A$ there exists $n\ge 1$ such that
$aA + \Ann_A (a^n) =A$,
then $A$ is finite-dimensional.
\end{prop}

\begin{proof}
It was shown in \cite{SSW} that $A$ is a  PI algebra
and its prime radical
${\mathfrak N} = {\mathfrak N}(A)$ is nilpotent
(moreover, it coincides with the Jacobson radical
$J=J(A)$, see, e.g., \cite{Ke}).
So $\bar A = A/{\mathfrak N}$ is a finitely generated
semiprime algebra which is left and right
 Goldie \cite{SW,SSW}.

Let $\bar a \in \bar A$ be a right regular element, i.e.,
$\bar a  \bar x = 0$  implies $\bar x=0$ in $\bar A$
(such an element necessarily exists in
any semiprime Goldie algebra).
Then $\Ann_{A} (a^n) \subseteq {\mathfrak N}$
for any $n\ge 1$, and
\[
aA + {\mathfrak N} = A.
\]
It is easy to derive that $\bar A$ contains a unit
and any (right) regular element is invertible.
Hence, $\bar A$ coincides with its classical quotient
algebra $Q(\bar A)$.
By the Goldie theorem, $\bar A=Q(\bar A)$ is semisimple Artinian,
so it is equal to a finite direct sum of
simple (and finitely generated) algebras of at most linear growth.
The main result of \cite{SW} implies that any algebra of this
kind is finite-dimensional, so
$\dim \bar A < \infty$.
 Finally, one may use
the Kaplansky theorem  and nilpotency of $\mathfrak N$
to
deduce that $A$ is finite-dimensional itself.
\end{proof}

The following statement generalizes the result of
\cite{Re1} (c.f. Theorem~\ref{thm5.2}).

\begin{prop}\label{prop5.3}
If there exists an element $a\in C$ such that
$a\oo{0} C = C$,
then $C\simeq \Cend_{N}$.
\end{prop}

\begin{proof}
First, consider
the homomorphism of left $C$-modules
$R: C\to C$, $x\mapsto \{x\oo{0} a\}$.
If $x\in \Kerr_C R$, then by (\ref{eq2.2.4}) we have
\[
0=\{x\oo{0} a\}\oo{\omega } C = x\oo{\omega } (a\oo{0} C) = x\oo{\omega }C,
\]
so $x=0$. Proposition~\ref{MainProp} and Remark \ref{rem4.1} imply
that either $C=\{C\oo{0} a\}$, or $C$ has a finite faithful
irreducible right module. In the last case, $C\simeq \Cend_{N,Q}$
by Theorem \ref{thmFF}.

By Lemma \ref{lem2.3.1}, $A_0=A_0(C)$ is finitely generated.
If the condition (\ref{eq5.1}) fails for some $b(0)\in A_0$,
then by Corollary \ref{cor4.1}
we have $C\simeq \Cend_{N,Q}$.
If (\ref{eq5.1}) holds for any $b(0)\in A_0$,
then $\dim A_0 <\infty$ by Proposition~\ref{prop5.2}.

Therefore, either $C$ is isomorphic to $\Cend_{N,Q}$
as in Theorem \ref{thmFF}, or $\{C\oo{0} a\} = C$ and $\dim A_0 <\infty$.
In the last case,
$\dim (C\oo{0}a)<\infty $, but $C=\{C\oo{0} a\}= C\oo{0}a  +DC$,
so $\dim C/DC < \infty $.
Then Lemma \ref{lem2.3.2} implies $C$ to be of finite type,
but $C$ is assumed to be of linear growth.

Hence, $C$ is necessarily isomorphic to $\Cend_{N,Q}$, $\det Q \ne 0$.
Since there exists $a\in C$ such that
$a\oo{0} C =C$, the matrix $Q$ has to be invertible, so
$\Cend_{N,Q}=\Cend_N$.
\end{proof}

\begin{lem}[c.f. \cite{Z2}]\label{lem5.3}
If there exists an element $a\in C$ such that
$a\oo{0}(a\oo{0}C) = a\oo{0}C \ne 0$,
then
$C_0 = a\oo{0} \{C\oo{0} a\}$
is a finitely generated conformal algebra
such that every proper ideal $I$ of $C_0$
satisfies $a\oo{0}I=0$.
\end{lem}

\begin{proof}
Since $C$ is simple, we have
\[
\{C\oo{0} a\} \oo{\omega } (a\oo{0} C) = C
\]
(the left hand side is an ideal which contains
$a\oo{0}a\oo{0}a\oo{0}a\ne 0$).
Any element $b\in C$ could be represented as
\begin{equation}\label{eq5.2}
b = \sum\limits_i \{x_i \oo{0} a\}\oo{n_i} (a\oo{0} y_i),
\quad x_i,y_i\in C.
\end{equation}
In particular,
every element $b\in B$ (where $B$ is a finite set of generators)
could be represented as (\ref{eq5.2}). So there exists a finite set $X\subset C$
such that
$B\subseteq \{X\oo{0} a\}\oo{\omega}(a\oo{0} X)$. Then the finite-dimensional
subspace
\[
a\oo{0}\{X\oo{0} a\} + a\oo{0} \{(X\oo{\omega} X) \oo{0} a\}
\]
generates $C_0$ as a conformal algebra.

Let $I$ be an ideal of $C_0$. Then
\[
J = \{C\oo{0} a\} \oo{\omega } I \oo{\omega } (a\oo{0} C)
\]
is an ideal of $C$, so either $J=0$ or $J=C$.
Note that
\[
a\oo{0} \{J \oo{0} a\}
= (a\oo{0} \{C\oo{0} a\})
 \oo{\omega } I \oo{\omega}
    (a\oo{0}\{C\oo{0} a\}) \subseteq I.
\]
Hence, $J=C$ implies $I=C_0$. If $J=0$, then
$ 0= \{C\oo{0} a\} \oo{\omega } I \oo{\omega } (a\oo{0} C)
    =C\oo{\omega} (a\oo{0} I)\oo{\omega } (a\oo{0} C)$
and
\begin{equation}\label{nilpI}
(a\oo{0} I)\oo{\omega } (a\oo{0} C) = 0.
\end{equation}

Now, consider an arbitrary element $x\in I$, $x=a\oo{0} \{ y \oo{0} a \}$.
For any $b\in C$ there exists an element $c\in C$
such that
$a\oo{0}b = a\oo{0}a\oo{0} c$.
Thus,
$
a\oo{0}x\oo{\omega} b = a\oo{0}a\oo{0} y \oo{\omega } (a\oo{0} b) =
a\oo{0}a\oo{0} (y \oo{\omega } (a\oo{0} a\oo{0} c)) =
a\oo{0}a\oo{0} (\{y \oo{0} a\} \oo{\omega} (a\oo{} c) ) =
(a\oo{0}x)\oo{\omega} (a\oo{0} c) =0
$
by (\ref{nilpI}).
Hence,
$a\oo{0} I $
annihilates the whole $C$, so $a\oo{0} I = 0$.
\end{proof}

We would like to state a simple proof of
Theorem \ref{thm5.3}, which is important for
further considerations.

\begin{proof}[Proof of Theorem \ref{thm5.3}]
It is sufficient to show that $C\simeq \Cend_{N,Q}$,
$N\ge 1$, $\det Q \ne 0$.
It is easy to note that $\Cend_{N,Q}$ contains an idempotent
if and only if the canonical diagonal form of $Q$
is of the form  $\diag(1,f_2,\dots, f_N)$.

Fix a finite set of generators $B\subset C$.
Let $e\in C$ be an idempotent.
Since
$e\oo{0}e = \{e\oo{0} e\} = e$,
the subalgebra
  $C_0 = e\oo{0} \{C \oo{0} e \}$ is simple and finitely generated by
 Lemma \ref{lem5.3}.
Moreover,
$C_0= H(e\oo{\omega } C \oo{\omega } e)$,
and $C_0$ contains a unit, e.g., $e = e\oo{0} \{ e\oo{0} e\}\in C_0$.
Hence, either $C_0$ is finite,
or $C_0\simeq \Cend_N$
for some $N\ge 1$
 (see Theorem \ref{thm5.2} or Proposition \ref{prop5.3}).

If $C_0$ is finite, then $C$ is finite itself.
Indeed, one may present
\[
C = HB + H(B\oo{\omega} B) + H(B\oo{\omega } C \oo{\omega } B)
\]
and note that there exists a finite-dimensional subspace $X \subset C$
such that
$B \subset X\oo{\omega } e \oo{\omega } X$
 (since $C\oo{\omega } e\oo{\omega} C$
 is a non-zero ideal of $C$).
Thus,
\[
C = HB + H(B\oo{\omega} B) + H(X \oo{\omega} C_0 \oo{\omega}  X)
\]
is a finite conformal algebra.

If $C_0\simeq \Cend_N$, then there exists an element $v\in C_0$ corresponding to
$x\Id_N \in M_N(\Bbbk[D,x])\simeq \Cend_N$. We may assume that there exists $n\ge 1$
such that (\ref{eq5.0}) holds for $v$, i.e.,
\[
v\oo{0} C  + \Kerr_C v(0)^n = C
\]
(otherwise, Corollary~\ref{cor4.1} implies the claim).
In particular, for the element
$a = v^{\oo{0} n}\equiv v\oo{0}\dots \oo{0}v \in C_0$
 (corresponding to $x^n\Id_N $)
there exists $y\in C$ such that
\begin{equation}\label{contr}
a\oo{0}a\oo{0} y = a\oo{0} e_1 = a.
\end{equation}
Here $e_1$ is the element of $C_0$ corresponding to $\Id_N \in \Cend_N$
(in fact, the isomorphism between $C_0$ and $\Cend_N$ could be chosen in such a way
that $e=e_1$, see \cite{Re1}).
Consider $z= e_1\oo{0} \{y\oo{0} e_1\}\in C_0$,
and note that
\begin{equation}\label{contr2}
a\oo{0}a\oo{0} z  = a\oo{0}a\oo{0}e_1\oo{0} \{y\oo{0} e_1 \}
= \{a\oo{0} e_1 \} = a.
\end{equation}
This relation is clearly impossible in $\Cend_N$.
\end{proof}

The following lemma is a particular case of the general statement from
\cite{Z2}.

\begin{lem}\label{Lifting}
Let $C$ be a conformal algebra,
and let $I$ be an ideal of $C$
such that $C\oo{\omega } I = 0$.
If $\overline{C} = C/I $ contains an idempotent
$\bar e$, then there exists a preimage
$e\in C$ which is also an idempotent.
\end{lem}

\begin{proof}
Let $e\in C$ be an arbitrary preimage of
the idempotent $\bar e \in \overline{C}$.
Then $e - e\oo{0}e \in I$,
and $e\oo{n} e \in I$ for any $n\ge 1$.
In particular, $e_1 = e\oo{0} e $
is also a preimage of~$\bar e$.
Since $C\oo{\omega } I =0$, we have
$e_1\oo{0} e_1 = e_1$,
and
$e_1\oo{n} e_1 = 0$
for any $n\ge 1$.
Therefore, $e_1$ is an idempotent of~$C$.
\end{proof}

\begin{proof}[Proof of Theorem \ref{thm5.1}]
Let us suppose that $C\not\simeq \Cend_{N,Q}$. Then  by
Corollary~\ref{cor4.1} for every $L\in {\mathcal L}(C)$ there exists
$n\ge 1$ such that $L(C) + \Kerr_C L^n = C$.
In particular, every element of $C$ satisfies the condition (\ref{eq5.0}).

Assume that there exists an element
$e\in C$ such that $e$ is not nilpotent
with respect to the $0$-product.
Then by (\ref{eq5.0}) there exists an integer $n\ge 1$ such that
\[
e\oo{0} C + \Kerr_C e(0)^n = C,
\]
so
\[
e(0)^n C = e(0)^n (e\oo{0} C) = e(0)^{n+1} C \ne 0.
\]
For $a=e(0)^{n-1}e \in C$ we have $a\oo{0} (a \oo{0}C) = a\oo{0} C\ne 0$.
By Lemma \ref{lem5.3}, the conformal algebra
$C_0 = a\oo{0} \{C\oo{0} a\}$
 is finitely generated, and for every proper ideal
$I$
of $C_0$ we have
$C_0\oo{\omega } I = 0$.
Moreover, there exists $x=a\oo{0} \{a\oo{0}a\}\in C_0$
such that $x\oo{0} C_0 = C_0$.

Let us choose a maximal ideal $I$ of $C_0$, and consider
the simple finitely generated associative algebra
$C_1 = C_0/I$. Since $\GKdim C_1\le \GKdim C =1$,
we may use Proposition \ref{prop5.3} and Theorem \ref{thmFF}
to conclude that $C_1$ is isomorphic either to $\Cend_N$ or
to $\Curr M_N(\Bbbk)$. In any case, $C_1$ contains a unit,
so $C_0$ has an idempotent
by Lemma~\ref{Lifting}.
Theorem \ref{thm5.3} implies that $C\simeq \Cend_{N,Q}$, $N\ge 1$,
$Q=\diag(1, f_2,\dots, f_N)$.
But this algebra clearly does not satisfy the condition (\ref{eq5.0}).

It remains to consider the case when every element of $C$ is nilpotent
with respect to the $0$-product. Let us choose any irreducible module $V$ of
$C$. As it was shown in \cite{BKL1} (see Proposition \ref{prop2.1}), there exist
$0\ne u\in V$ and $\gamma \in \Bbbk $ such that $\{C {}_{\gamma} u\} =V$.
In particular, there exists an element $b\in C$ such that
$\{b {}_{\gamma} u\} = u$.
Note that
\[
u= \{b {}_{\gamma} \{b {}_{\gamma} \dots {}_{\gamma} \{b_{\gamma} u \}\dots \}\}
= \{\{b {}_{0}\dots {}_{0} b\} _{\gamma} u \} =0
\]
(see (\ref{eq2.2.2}) and (\ref{eq2.2.3})),
so $u=0$ in contradiction with $V\ne 0$.
\end{proof}

\begin{thm}[\cite{BKL1}]\label{thmUNQ}
Two conformal algebras $\Cend_{N_1,P}$ and $\Cend_{N_2,Q}$
$(\det P(x),\break \det Q(x)\ne 0)$
 are  isomorphic if and only if
$N_1=N_2=N$  and there exists $\alpha \in \Bbbk $ such that $P(x)$
and $Q(x+\alpha) $ have the same canonical diagonal form, i.e.,
there exist matrices $A,B\in M_N(\Bbbk[x]) $,
$\det A(x),\det B(x)\in \Bbbk\setminus \{0\}$, such that
$P(x)A(x) = B(x)Q(x+\alpha)$.
\qed
\end{thm}

Theorems \ref{thm5.1} and \ref{thmUNQ} provide the complete description
of simple finitely generated associative conformal algebras of
linear growth up to isomorphism.

\section*{Acknowledgement}
I am grateful to L.~Bokut,  E.~Zel'manov, Yu.~Mal'cev for
helpful discussions and consultations. I would like to acknowledge
the referee for many remarks and suggestions that led to improvement
of the presentation and language.

\end{document}